\documentclass[a4paper,12pt,leqno]{amsart}   

\usepackage{amsmath, amssymb,latexsym}
\usepackage{amsthm,amscd,amsfonts}
\usepackage{graphicx,psfrag} 

\usepackage[english]{babel}
\usepackage[latin1]{inputenc}

\usepackage{verbatim}
\usepackage[all, cmtip]{xy}
\usepackage{extpfeil}
\usepackage{bbm}
\usepackage{paralist}
\usepackage[T1]{fontenc}
\usepackage[colorlinks=false,citecolor=magenta,pagebackref=true,urlcolor=magenta,pdftex]{hyperref}
\usepackage{nicefrac}
\usepackage{microtype}
\usepackage{epstopdf}

\usepackage{mathrsfs}

\newcommand{\cm}{Cohen-Macaulay}
\newcommand{\cmc}{Cohen-Macaulay complex}

\begin{document}




\newtheorem{thm}{Theorem}[section]
\newtheorem{prop}[thm]{Proposition}
\newtheorem{lem}[thm]{Lemma}
\newtheorem{cor}[thm]{Corollary}
\newtheorem{conj}[thm]{Conjecture}
\newtheorem{ddef}[thm]{Definition}
\newtheorem{ex}[thm]{Example}
\newtheorem{rem}[thm]{Remark}
\newtheorem{notation}[thm]{Notation}

\numberwithin{equation}{section}
\numberwithin{figure}{section}

\newcommand{\bthm}{\begin{thm}} \newcommand{\ethm}{\end{thm}}
\newcommand{\bthms}{\begin{thm*}} \newcommand{\ethms}{\end{thm*}}
\newcommand{\blem}{\begin{lem}} \newcommand{\elem}{\end{lem}}
\newcommand{\bcor}{\begin{cor}} \newcommand{\ecor}{\end{cor}}
\newcommand{\bprop}{\begin{prop}} \newcommand{\eprop}{\end{prop}}
\newcommand{\bproof}{\begin{proof}} \newcommand{\eproof}{\end{proof}}
\newcommand{\bddef}{\begin{ddef}} \newcommand{\eddef}{\end{ddef}}
\newcommand{\bconj}{\begin{conj}} \newcommand{\econj}{\end{conj}}
\newcommand{\brem}{\begin{rem}} \newcommand{\erem}{\end{rem}}
\newcommand{\bca}{\begin{cases}} \newcommand{\eca}{\end{cases}}

\newcommand{\att}[1]{\tcr{{\em #1}}}

\newcommand{\beq}{\begin{equation}} \newcommand{\eeq}{\end{equation}}
\newcommand{\beqs}{\begin{equation*}} \newcommand{\eeqs}{\end{equation*}}
\newcommand{\beqa}{\begin{eqnarray}} \newcommand{\eeqa}{\end{eqnarray}}
\newcommand{\beqas}{\begin{eqnarray*}} \newcommand{\eeqas}{\end{eqnarray*}}
\newcommand{\barr}{\begin{array}} \newcommand{\earr}{\end{array}}
\newcommand{\btab}{\begin{tabular}} \newcommand{\etab}{\end{tabular}}

\newcommand{\bit}{\begin{itemize}} \newcommand{\eit}{\end{itemize}}
\newcommand{\ben}{\begin{enumerate}} \newcommand{\een}{\end{enumerate}}
\newcommand{\bce}{\begin{center}} \newcommand{\ece}{\end{center}}

\newcommand{\defeq}{\stackrel{\rm def}{=}}
\newcommand{\noset}{\varnothing}
\newcommand{\bd}{\partial}

\newcommand{\ho}{\hat{1}}
\newcommand{\hz}{\hat{0}}
\newcommand{\whz}{\widehat{0}}
\newcommand{\who}{\widehat{1}}
\newcommand{\wh}{\widehat}
\newcommand{\wt}{\widetilde}
\newcommand{\cover}{\prec}
\newcommand{\cov}{\lhd}
\newcommand{\cocov}{\rhd}
\newcommand{\opp}{\mathrm{opp}}
\newcommand{\join}{\vee}
\newcommand{\meet}{\wedge}

 
\renewcommand{\AA}{\mathcal{A}}
\newcommand{\BB}{\mathcal{B}}
\newcommand{\CC}{\mathcal{C}}
\newcommand{\DD}{\mathcal{D}}
\newcommand{\EE}{\mathcal{E}}
\newcommand{\FF}{\mathcal{F}}
\newcommand{\HH}{\mathcal{H}}
\newcommand{\II}{\mathcal{I}}
\newcommand{\LL}{\mathcal{L}}
\newcommand{\MM}{\mathcal{M}}
\newcommand{\NN}{\mathcal{N}}
\newcommand{\PP}{{\mathcal{P}}}
\newcommand{\NP}{\mathcal{NP}}
\newcommand{\QQ}{\mathcal{Q}}
\newcommand{\RR}{\mathcal{R}}
\renewcommand{\SS}{\mathcal{S}}


\newcommand{\C}{\mathbb{C}}
\newcommand{\F}{\mathbb{F}}
\renewcommand{\L}{\mathbb{L}}
\newcommand{\N}{\mathbb{N}}
\renewcommand{\P}{\mathbb{P}}
\newcommand{\Q}{\mathbb{Q}}
\newcommand{\R}{\mathbb{R}}
\newcommand{\Z}{\mathbb{Z}}

\newcommand{\kk}{\mathbf{k}}
\newcommand{\tH}{\wt{H}}
\newcommand{\lk}{\mathrm{link}}


\newcommand{\al}{\alpha}
\newcommand{\be}{\beta}
\newcommand{\de}{\delta} \newcommand{\De}{\Delta} 
\newcommand{\ga}{\gamma} \newcommand{\Ga}{\Gamma}
\newcommand{\la}{\lambda} \newcommand{\La}{\Lambda}
\newcommand{\om}{\omega} \newcommand{\Om}{\Omega}
\newcommand{\si}{\sigma} \newcommand{\Si}{\Sigma}
\newcommand{\ze}{\zeta}
\newcommand{\vphi}{\varphi}
\newcommand{\eps}{\varepsilon}

\newcommand{\X}{\bf{X}}
\newcommand{\norm}[1]{\lVert#1\rVert}
\newcommand{\abs}[1]{\lvert#1\rvert}
\newcommand{\qbinom}[2]{\left[\ba{c}{#1}\\{#2}\ea\right]} 
\newcommand{\comp}{\models}

\newcommand{\vanish}[1]{}

\newcommand{\st}{\,:\,} 
\newcommand{\sbseq}{\subseteq}
\newcommand{\spseq}{\supseteq}


\newcommand{\larr}{\leftarrow}
\newcommand{\rarr}{\rightarrow}
\newcommand{\Larr}{\Leftarrow}
\newcommand{\Rarr}{\Rightarrow}
\newcommand{\lrarr}{\leftrightarrow}
\newcommand{\Lrarr}{\Leftrightarrow}

\newcommand{\longlarr}{\longleftarrow}
\newcommand{\longrarr}{\longrightarrow}
\newcommand{\Longlarr}{\Longleftarrow}
\newcommand{\Longrarr}{\Longrightarrow}
\newcommand{\longlrarr}{\longleftrightarrow}
\newcommand{\Longlrarr}{\Longleftrightarrow}

\title[]{``Let $\De$ be a Cohen-Macaulay complex $\ldots$'' }
\author[Anders Bj\"orner]{Anders Bj\"orner}

\dedicatory{\small Dedicated to Richard Stanley on the occasion of his 70th birthday\\
\vspace*{-27mm}
\begin{center}
\includegraphics[scale=0.5]{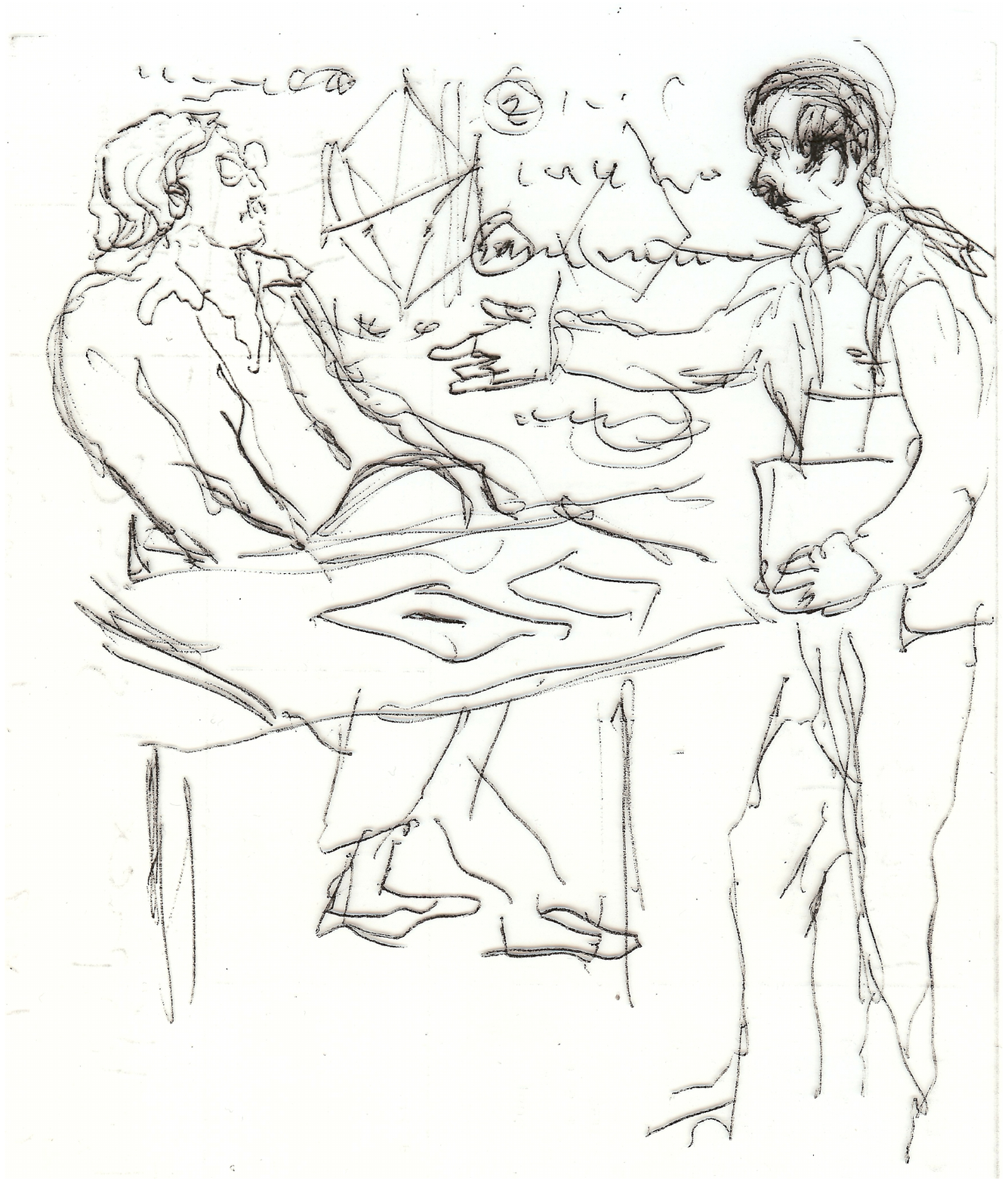} 
\end{center}
\vspace*{-33mm}
}

\address{Kungl. Tekniska H\"ogskolan, Matematiska Inst.,
 S-100 44 Stockholm, Sweden}
\email{bjorner@math.kth.se}
\thanks{Based on a talk given at the ``Stanley@70'' birthday conference'', MIT, June 2014}


\begin{abstract} 

The concept of Cohen-Macaulay complexes emerged in the mid-1970s
and swiftly became the focal point of
an attractive and richly connected new area of mathematics,
at the crossroads of combinatoics, commutative algebra and topology.
As the main architect of these developments, Richard Stanley has made fundamental  contributions
over many years. 


This paper contains some brief mathematical discussions
related to the Cohen-Macaulay property, and some personal memories.
The characterization of Gorenstein*
and homotopy Gorenstein* complexes and the relevance  in that connection of the Poincar\'e conjecture is discussed.
Another topic is combinatorial aspects of a recent result on
the homotopy Cohen-Macaulayness of certain subsets of 
geometric lattices, motivated by questions in tropical geometry.
\end{abstract}

\maketitle

\section{Introduction}

As is well known, Richard has a very good sense of humor. One of his favorite
jokes around 1980 went somewhat like this: 
\begin{quote}
{\em ``Can you imagine that at some time in the future
it will be possible to begin a general math lecture with 
{``Let $\De$ be a Cohen-Macaulay complex $\ldots$''}
and go on from there without further explanation.''}
\end{quote}

The fact that this would not work as a joke today 
testifies to 
the emergence of this concept from
what was then a remote corner of combinatorics into
mainstream mathematics.


The concept of Cohen-Macaulay complexes arose in the mid-1970s from the work of
Melvin Hochster \cite{Ho72,Ho77}, Gerald Reisner \cite{Re75}, 
Richard Stanley \cite{St75,St77}, and (for the case of posets) 
Kenneth Baclawski \cite{Ba76,Ba80}.
Others soon followed and an attractive new area of mathematics took shape. 
The theory of \cm\ complexes 
has applications to several concepts in combinatorics,
such as matroids, polytope boundaries, geometric lattices, buildings, intersection 
lattices of hyperplane arrangements, and more.

Much of the appeal of the concept of CM-ness stems from its ``interdisciplinary'' character, building on and having bearing on several mathematical areas in addition to combinatorics, notably several parts of algebra, geometry and topology. However, what bestows lasting  importance on the concept   is its remarkable record of being a key ingredient both for ``abstract'' theoretical understanding and for ``concrete'' problem solving in several of these diverse areas.

As the main architect of these developments, Richard has made fundamental  contributions over many years. His proof in 1975 of
the Upper Bound Conjecture for spheres  \cite{St75} catapulted Cohen-Macaulay 
complexes into the limelight. Then his 1977 paper \cite{St77}
outlined the contours of a theory with many beautiful
results and many appealing problems.
His contributions to this area 
have been essential for the development of  algebraic combinatorics, 
and have had a significant impact also on a wider mathematical  territory, 
particularly for ring theory. 




\medskip

In this paper I meander among basic facts, 
brief mathematical discussions, and personal memories.
One topic discussed is the topological characterization of homotopy Gorenstein* complexes and  the 
relevance of the Poincar\'e conjecture in that connection,
observed a few years ago by P. Hersh and the author (Section 4).
Another topic is a recent result of K. Adiprasito and the author on
the Cohen-Macaulayness of certain subsets of 
geometric lattices (Section 3), motivated by questions in tropical geometry.
This result can be illustrated\footnote{in the sense of Diaconis
\cite{Di14}.}  by the following small ``story''.

Suppose that each street corner of a crime-ridden city has 
been given a 
ranking number,
reflecting how safe it is to visit that corner. To normalize the grading,
the average rank has been set to be zero. Street corners with
positive rank are considered safe, those with negative rank are not.
Furthermore, a street in the city is considered
safe if the average rank of all street corners along that street is positive,
otherwise it is dangerous.
The question is: Is it possible to walk from any safe street corner
to any other safe one without ever passing a dangerous corner or walking along 
a dangerous street? For the answer, see Section \ref{CombAsp}.

\section{Cohen-Macaulay complexes revisited}



The purpose of this section is to remind the reader of some definitions and 
basic facts.

\subsection{Simplicial complexes everywhere}

We begin with a few words to  fix notation and agree on basic definitions.

A {\em simplicial complex} (or {\em abstract simplicial complex},
or just {\em complex})  is a finite set $V$  together
with a family $\De$ of subsets of $V$ such that $A \sbseq B\in \De$
implies that $A\in \De$. The elements of $V$ are the {\em vertices} and
the members of $\De$ the {\em faces} of the complex. 
We assume that the empty set is
a face.
The nonempty faces are the {\em proper} faces.
 The vertices are usually clear from context, and if so we use only
$\De$ to denote the complex.

The {\em dimension} of a face is one less that its cardinality, and the
dimension of $\De$ is the maximal dimension of any of its faces.

\medskip

A mathematical concept could hardly be simpler than that of a simplicial complex.
All that is required is a family of subsets of a finite set, closed under the
operation of taking subsets.
Simplicial complexes arise everywhere in combinatorics under
different names: hypergraphs, hereditary set families, order ideals in
the Boolean lattice, etc.

To define a Cohen-Macaulay complex is not as elementary.
There are two ways to proceed, equivalent as it turns out, 
via algebra (commutative rings) or via 
topology (simplicial homology).

\subsection{\cmc es via commutative rings}

 In commutative algebra simplicial complexes correspond to squarefree
 monomial ideals. In the following, {\bf k} denotes a field or the ring $\Z$ of integers.
 Suppose that  $\De$ is a complex on
 vertex set $V=\{1,2, \ldots, n\}$. Let $A=\kk[x_1, x_2, \ldots, x_n]$
 and let $I_{\De}$ be the ideal generated by squarefree monomials
 $x_{i_1}x_{i_2} \cdots x_{i_k}$ such that
$\{i_1 ,  i_2  , \ldots   i_k\} \notin\De$. The ideal $I_{\De}$ is called the
{\em Stanley-Reisner ideal} and the ring $\kk[\De] = A / I_{\De}$ the
{\em Stanley-Reisner ring}.

Let $$ 0 \rarr F_j  \rarr F_{j-1} \rarr\cdots\rarr F_0 \rarr A / I_{\De} \rarr 0$$
be a minimal free resolution of $\kk[\De]$ as an $A$-module. We
know from Hilbert's syzygy theorem that its length $j$ satisfies
$j\le n$, and from the Auslander-Buchsbaum theorem that $n-d \le j$, where
$d\defeq \dim \De +1$. The integer $n-j$ is called the
{\em depth} of $\kk[\De]$, and the rank of the free module $F_j$ is the {\em type} of 
$\kk[\De]$. Clearly, depth  $\le d$.

Here are the basic definitions, for the ring $\kk[\De]$ and for the complex $\De$. 
\bddef (i) $\kk[\De]$ is {\em Cohen-Macaulay} if its depth equals $d$.\\
(ii) $\kk[\De]$ is {\em Gorenstein} if it is a Cohen-Macaulay ring of type $1$.\\
(iii) $\De$ is {\em Cohen-Macaulay}, resp. {\em Gorenstein}, if 
$\kk[\De]$ is. \\
(iv) $\De$ is {\em Gorenstein*} if it is Gorenstein and not  acyclic over $\kk$.
\eddef

In the sequel we consider only
Gorenstein* complexes, the reason being that a general Gorenstein complex is a 
multiple cone over a  Gorenstein* complex, and the latter carries all 
relevant information.

\subsection{\cmc es via simplicial homology}

\medskip 
In topology simplicial complexes play the important role of encoding compact topological spaces. With a complex $\De$ is associated its 
{\em geometric realization}
 $||\De||$. Encoding reasonably nice compact spaces via triangulation 
 has always been an important tool in topology. This importance 
 has increased in recent years as triangulation is a necessary part of
 protocols for communicating with computers about topological spaces.

If $F$ is a face of $\De$, then 
$$\lk_{\De}(F) \defeq \{G\in \De \mid G\cup F\in \De \mbox{ and }
G \cap F =\emptyset\}
$$
is a subcomplex, called the {\em link of $\De$ at $F$}.  It carries local information 
about the complex
and its geometric realization. 
Note that 
$\lk_{\De}(\emptyset) =\De$. We say that $\lk_{\De}(F)$ is {\em proper} if $F\neq \emptyset$.

The following topological characterizations of 
\cm ness, 
due to Reisner \cite{Re75}, is a corner stone of the theory.

\bthm \label{Thm.Re75}
A simplicial complex $\De$ is {\em \cm\ over} {\bf k}
$$ \Longleftrightarrow \hspace{2mm}  \tH _i (\lk_{\De} (F) ; {\kk})=0 \mbox{ for all $F\in\De$ and all $ i < \dim(\lk_{\De}\ F)$}. $$
\ethm

\subsection{Commutative Algebra $\lrarr$ Algebraic Topology}

\vspace{7mm}
\begin{figure}[h]
\begin{center}
\includegraphics[scale=0.3]{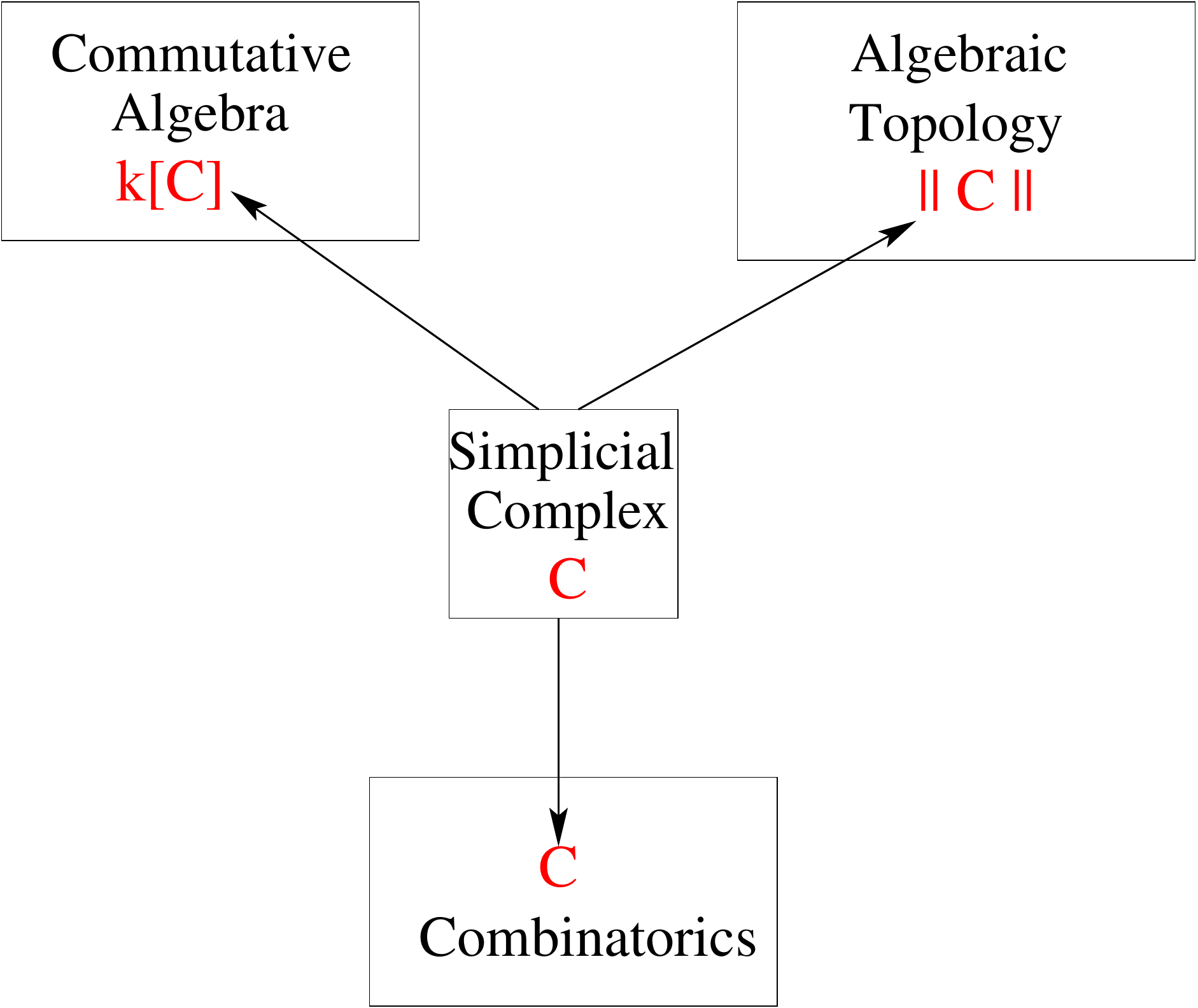} 
\caption{Simplicial complexes provide the link}
\label{fig11}
\end{center}
\vspace{7mm}
\end{figure}

At the core of Stanley-Reisner theory stands the discovery,
well exemplified by Reisner's Theorem \cite{Re75},
that some  central properties of commutative rings correspond
to important topological properties of simplicial complexes. 
 The following  formula,
due to Hochster \cite{Ho77},
gives a  particularly beautiful instance of this correspondence.

\label{hochster}
\beq  \be_{i,j} ( \mathbf{k} [\De]) = \sum_{E\sbseq V :  |E|=j} \dim_{\mathbf{k}}
\wt{H}_{j-i-1} (\De |_{E} : \mathbf{k})
\eeq

Here on the left hand side are the doubly indexed 
ring-theoretic Betti numbers,
which record the dimensions of the resolvants in a minimal free resolution
together with information about shifting of degrees. On the right hand
side are the topological Betti numbers of simplicial reduced homology of 
subcomplexes of $\De$ induced on subsets $E$ of the vertices.

A very interesting fact is that the ring-theoretic properties 
of being \cm\ or Gorenstein*  over $\mathbf{k}$ are
topological properties; they depend only on the homeomorphism type of
the space   $||\De ||$.
On the other hand, as pointed out by Reisner \cite{Re75} for the \cm\ case, these 
properties depend on field characteristic.
For instance, triangulations of the real projective plane $\R P^2$
are \cm\ over a field $\mathbf{k}$ if and only if $\mathrm{char}(\mathbf{k})\neq 2$.
Similarly, triangulations of real projective $3$-space  $\R P^3$
are Gorenstein* if $\mathrm{char}(\mathbf{k})\neq 2$,
but are not even \cm\ if $\mathbf{k}=\Z$ or $\mathbf{k}$ is a field of characteristic $2$.




\subsection{The homotopy \cm\ property }

We now  strengthen the concept of \cm ness by replacing the vanishing
of homology groups in the definition  by the vanishing of corresponding homotopy 
groups $\pi_i$.

\bddef
A simplicial complex $\De$ is {\em homotopy  \cm\ } if \\
$\pi _i (\lk_{\De}\ F)=0$ for all faces $F$ and all $ i < \dim(\lk_{\De}\ F)$.
\eddef

By basic algebraic topology we have the following characterization.
\bprop \label{Thm.HCM}
A simplicial complex $\De$ is {\em homotopy \cm}
$$ \Longleftrightarrow \hspace{2mm} \bca \De \mbox{ is \cm\ over $\Z$, and } \\
\mbox{all links of dimension $\ge 2$ are simply connected.} 
\eca
$$
\eprop

The concept of homotopy Cohen-Macaulay\-ness 
first apeared in Quillen's paper \cite{Qu78} on $p$-subgroups.
He offered this motivation: 
``we use the stronger  definition so that our theorems are in their best form''.

In spite of being extremely natural from a topological
point of view, this property is itself
not topologically invariant. This is witnessed by the $5$-dimensional sphere,
which in addition to its standard  triangulations, e.g. as the boundary of 
 a $6$-simplex, admits triangulations which are not 
homotopy \cm\ \cite{Ed75}. 

\begin{figure}[h]
\begin{center}
\includegraphics[scale=0.4]{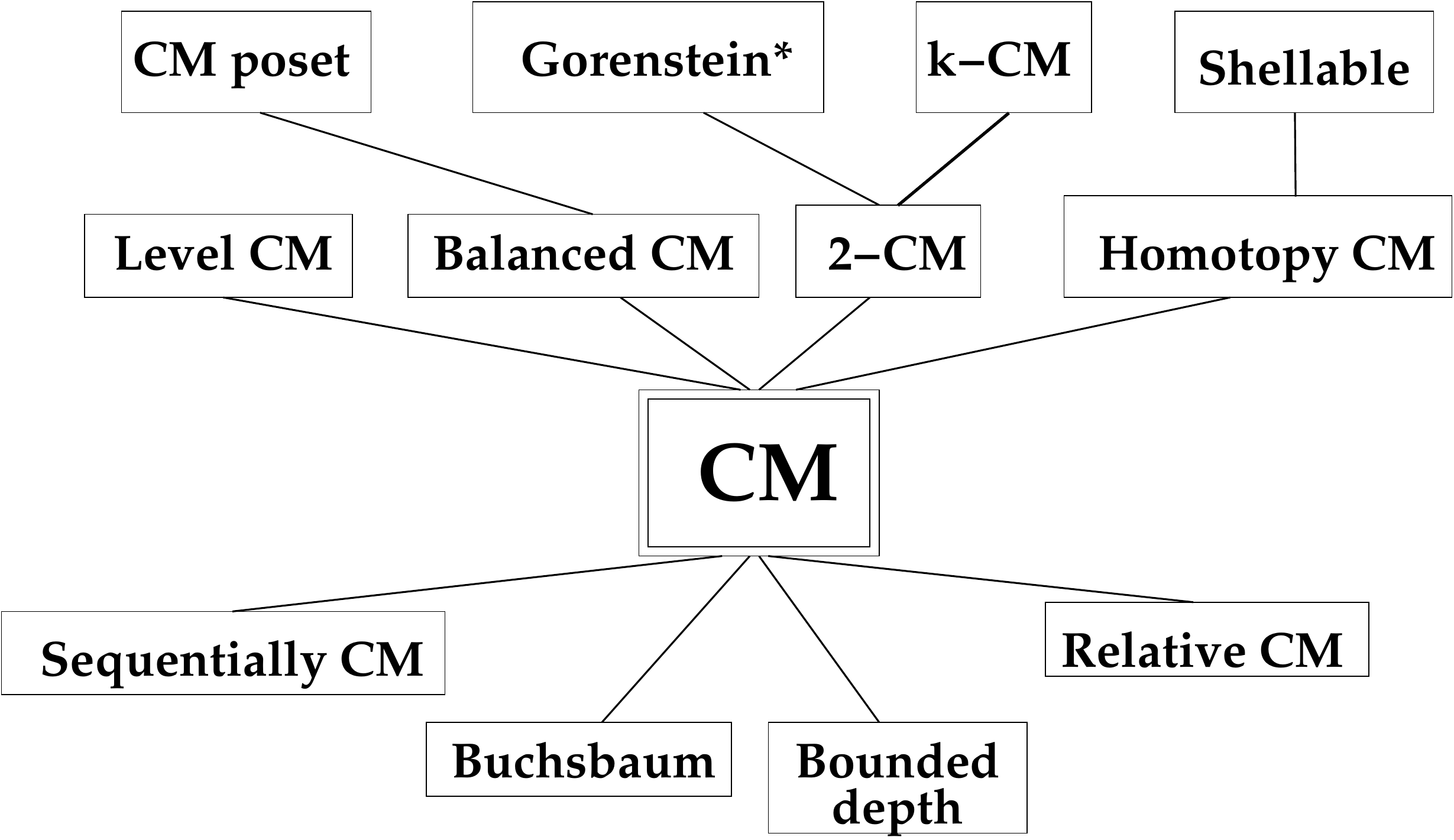}  
\caption{Properties of complexes  related to CMness.}
\label{fig1}
\end{center}
\vspace{7mm}
\end{figure}

Over the years many properties related to Cohen-Macaulayness
have been introduced and studied.
Figure \ref{fig1} shows some of the most important such properties, ordered
as a poset with logical implication arrows directed down.

This is not the place to enter a technical discussion
or attempt any kind of survey. 
As was stated initially, all that this section wants to convey are some brief reminders and 
comments about the basics.

\section{Cohen-Macaulayness of  filtered geometric lattices}


In this section we discuss a recent result motivated by questions in tropical geometry.
It concerns the Cohen-Macaulayness of $\De(P)$ for certain subsets $P$ of geometric  lattices.
Here the simplicial complex $\De(P)$  associated to a poset $P$
is the collection of its chains (totally ordered subsets). We do not distinguish 
notationally between a poset $P$ and its order complex $\De(P)$.

One of the properties of 
\cm ness that makes it so useful is its {\em resilience} -- several useful constructions
on complexes preserve the Cohen-Macaulay property. A good example is 
so called ``rank-selection'', by which is meant the removal of 
all elements of specified rank-levels in a \cm\ poset.
This operation preserves Cohen-Macaulayness, as was shown in varying degrees of
generality by Baclawski, Munkres, Stanley and Walker. 

By comparison, the class of Gorenstein* complexes is not at all resilient.
No removal of vertices leads to a Gorenstein* complex
of the same dimension.


\subsection{Geometric lattices}


A finite lattice is called {\em geometric} if it is
{\em semimodular} and {\em atomistic}. A geometric
lattice is pure, so it has a well-defined rank function,
$\rho: L\rarr \Z^+$. Being semimodular means that 
$$\rho(x) + \rho(y) \ge  \rho(x\wedge y) + \rho(x\vee y )$$
for alll $x,y\in L$,
and being atomistic means that every element in
the lattice is a join of atoms (atoms being the
elements that cover $\hz$).
Geometric lattices are  cryptomorphic to matroids, they
arise as the lattice of closed sets of a matroid, see \cite{Ox11}.

The order complex of a geometric lattice was one of the first examples 
of a  \cm\ complex. This came about in the following way. 
In his influential paper on the combinatorial M\"obius function \cite{Ro64},
Rota proved that the coefficients of the characteristic polynomial 
$$\chi(\LL; z)\defeq \sum_{x\in\LL}^{} \mu(\hat{0}, x)z^{r - \rho(x)}$$
of a geometric lattice $\LL$ alternate in sign:
$$ \chi(\LL; z) = z^r -a_{r-1} z^{r-1} + a_{r-2} z^{r-2}- a_{r-3} z^{r-3}+ a_{r-4} z^{r-4}- \cdots$$
 For the special case of chromatic polynomials of graphs
this was a well-known phenomenon.
The sign alternation is implied by the stronger set of inequalities
$$(-1)^{\rho(x)} \mu_{\LL}(\hat{0}, x) > 0 $$
valid for each $x\in \LL$.

 Rota's paper, which
deals with the Euler characteristic of order complexes, hinted at
a homological explanation of the result.
This was worked out  in a follow-up paper by Folkman \cite{Fo66}, who explicitly 
determined the homology of a geometric lattice. He showed that there is 
non-vanishing homology only in the top dimension. Since intervals
in a geometric lattice are themselves geometric lattices, 
Cohen-Macaulayness can be deduced.

\subsection{Filtered geometric lattices}

The resilience of the \cm\  property is well exemplified on the
class of geometric latices. For instance, it is known since long that
Cohen-Macaulayness and dimension are preserved by removal of any chain, and of certain antichains, from a geometric lattice \cite{Bj80, Ba82},
 and also by removal of any principal filter \cite{WW86}.
 The following is a recent discovery in this vein.

Let $\LL$ be a geometric lattice of rank $r$ and with set of atoms $A$,
and let $\om : A\rarr \R$ be a real-valued function assigning a nummber, called 
a {\em weight}, to each atom. Extend the weight function to all $X\in \LL$
by summation: $\om(X)=\sum_{a\in X} \om(a)$.
We assume that $\om$ is {\em generic}, meaning that 
$\emptyset \neq X\neq Y\neq\emptyset \Rarr 
\om(X) \neq \om(Y)$.

For $t\in \R$, let $\LL^{> t} \defeq \{X\in \LL \mid \om(X)>t\}$
with the induced partial order.
We call posets of this form {\em filtered geometric lattices}. 
They do not need to be lattices. The following result was
proved by Adiprasito and the author \cite{AB14}.

\begin{thm}\label{filter}
Suppose that $t \le \min[0, \om{(A)}]$. Then,
$\LL^{> t}$  is homotopy \cm\ and of the same dimension.
\end{thm}

The theorem is proved by techniques from topological
poset theory such as lexicographic shellability and Quillen-type fiber
arguments. It
was conjectured by Michalkin and Ziegler  \cite{MZ08}
in a stronger form, namely that $\LL^{> t}$
itself is shellable. This claim is still open.
The reason for their conjecture is that it implies
information that is crucial for proving Lefschetz-type
section theorems in tropical geometry. 
We will not here pursue this path into the territory of tropical geometry.

From now on we specialize to the conditions  $t=\om(A)=0$. Then,
with respect to the weights  $\om$ and $-\om$,
we have that
$\LL$ splits into two parts:
$$ \LL^+ \defeq\LL^{>0}  \mbox{ \quad and \quad } \LL^- \defeq\LL^{<0}
$$
Both parts are homotopy \cm\ and of the same dimension as $\LL$.

All properties of geometric lattices that depend only on \cm ness
can now automatically  be generalized to filtered geometric lattices.
For example, we get filtered characteristic polynomials,
$$\chi(\LL^+; z)\defeq \sum_{x\in\LL^+}^{} \mu_{\LL^+}(\hat{0}, x)z^{r- \rho(x)}.$$
whose coefficients alternate in sign:
$$ \chi(\LL^+; z) = z^r -b_{r-1} z^{r-1} + b_{r-2} z^{r-2}- b_{r-3} z^{r-3}+ b_{r-4} z^{r-4}- \cdots$$

For instance, the weighted configuration in Figure \ref{fig3v}
has these filtered characteristic polynomials:
$$ \chi(\LL^+; z) = z^3 -4 z^{2} + 3 z \mbox{ and }
\chi(\LL^-; z) = z^3 -2 z^{2} +  z.
$$

The sign alternation is implied by 
$$(-1)^{\rho(x)} \mu_{\LL^+}(\hat{0}, x) \ge 0, $$
which is valid for each $x\in \LL$ as a consequence of
\cm ness.

As a special case,
putting generic weights on the edges of a graph we are
led to consider   ``filtered chromatic polynomials''.
Such polynomials have to our knowledge not been explored.

\subsection{Combinatorial aspects}\label{CombAsp}

We will now have a look at what Theorem \ref{filter} says
in the first non-trivial case; that of rank $3$. The proper part of a geometric lattice 
of rank $3$ is a bipartite graph, and for graphs, Cohen-Macaulayness just means
being connected. This connectivity result can be illustrated in the 
following concrete way for matroids representable over $\R$.

\begin{figure}[h]
\begin{center}
\includegraphics[scale=0.16]{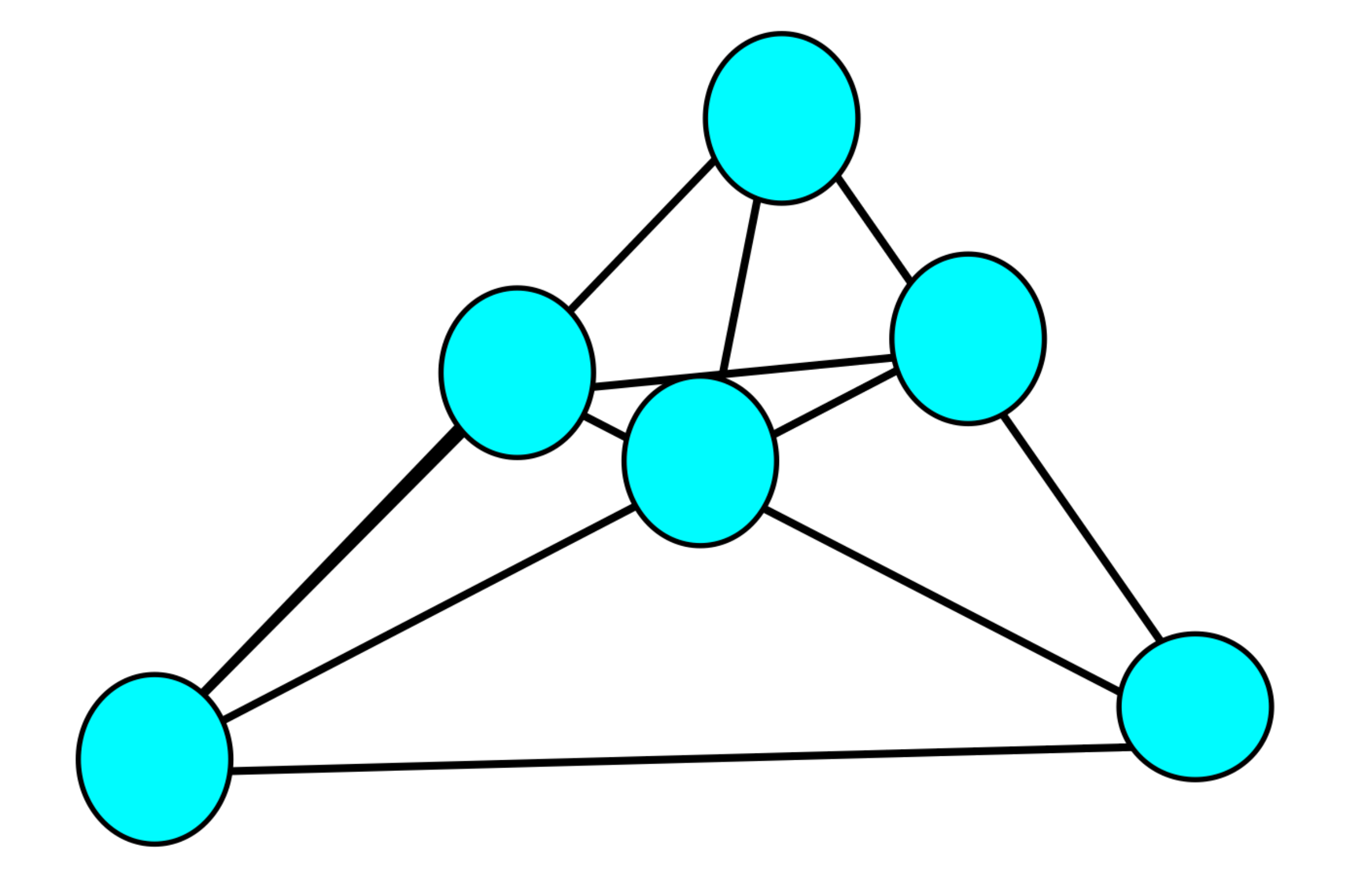}  
\vspace*{-2mm}
\caption{A configuration of $6$ points spanning $7$ lines in $\R^2$.}
\label{fig3x}
\end{center}
\vspace{2mm}
\end{figure}

Consider a finite collection of points in the plane $\R^2$. Along with the points we consider also the lines that are spanned by subsets of the points,
 see Figure \ref{fig3x}.
 
\begin{figure}[h]
\begin{center}
\includegraphics[scale=0.2]{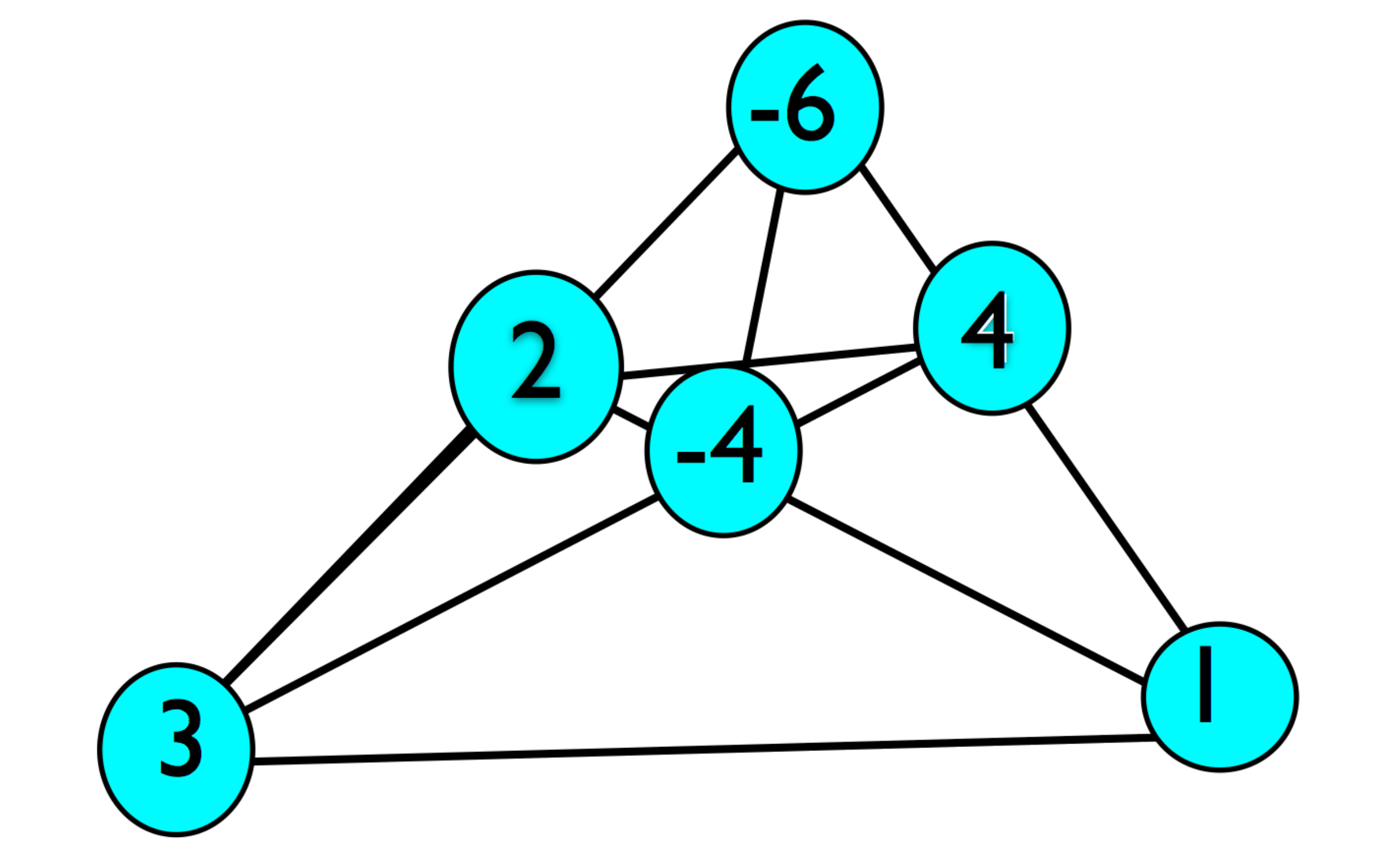}  
\vspace*{-2mm}
\caption{The weighted point configuration.}
\label{fig3v}
\end{center}
\vspace{3mm}
\end{figure}

A real number, called its {\em weight}, is assigned to each point, and we assume that
these numbers sum to zero, but are otherwise generic  (Figure \ref{fig3v}).

\begin{figure}[h]
\begin{center}
\includegraphics[scale=0.2]{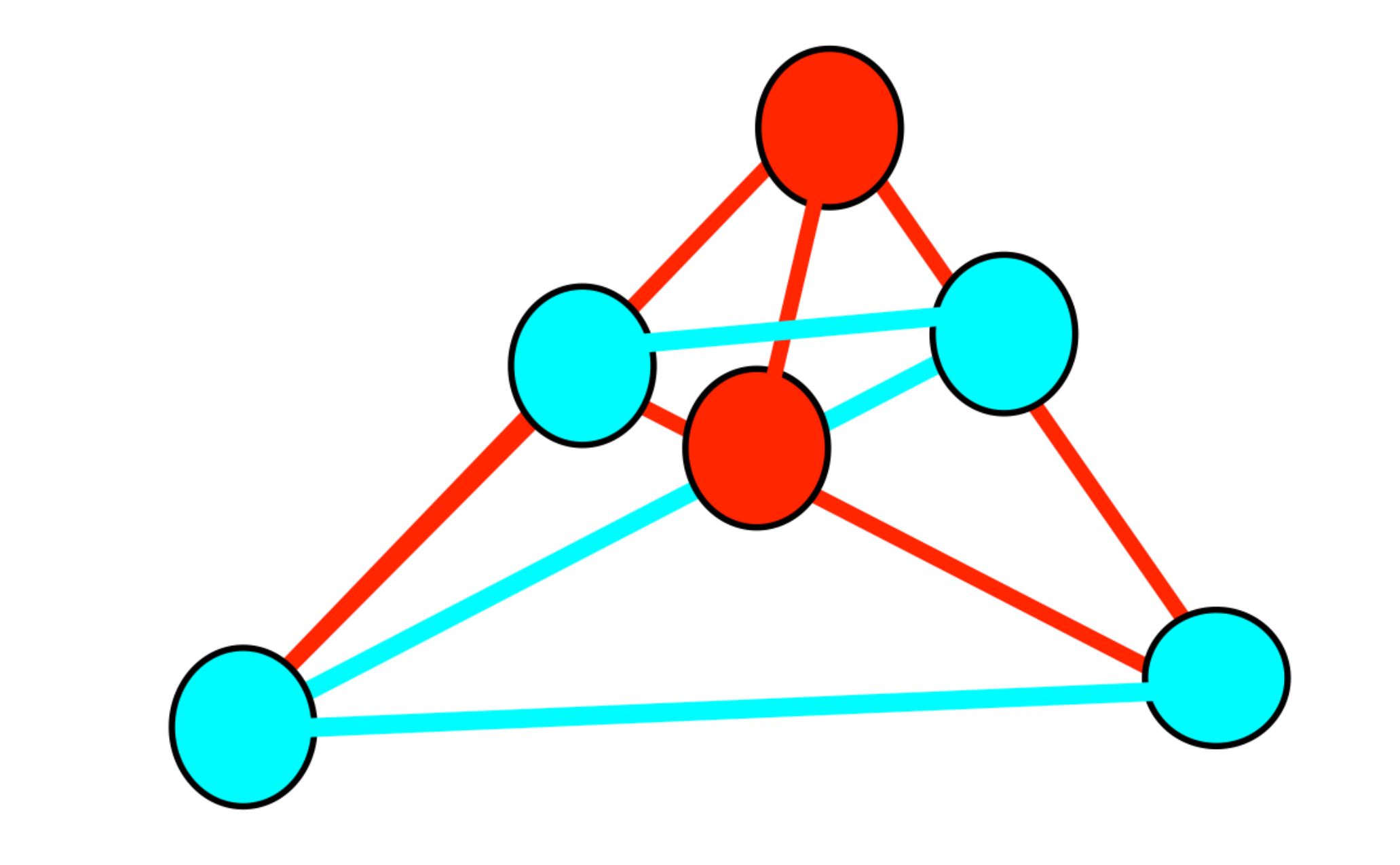}  
\vspace*{-5mm}
\caption{Negative points and lines are marked red.}
\label{fig3y}
\end{center}
\vspace{3mm}
\end{figure}

By the weight of a line we mean the sum of the weights of its points.
With this  we can distinguish between positive and negative lines.
For instance, in our running example we 
have four negative lines  $\{1, -4, 2\}$, $\{1, 4, -6\}$, $\{2, 3, -6\}$, $\{-4, -6\}$.
This leaves three positive lines (Figure \ref{fig3y}).


The question now is this: Is it possible to walk in the configuration from any positive point 
to any other positive point, never visiting a negative point or walking
along a negative line? 

\begin{figure}[h]
\begin{center}
\includegraphics[scale=0.2]{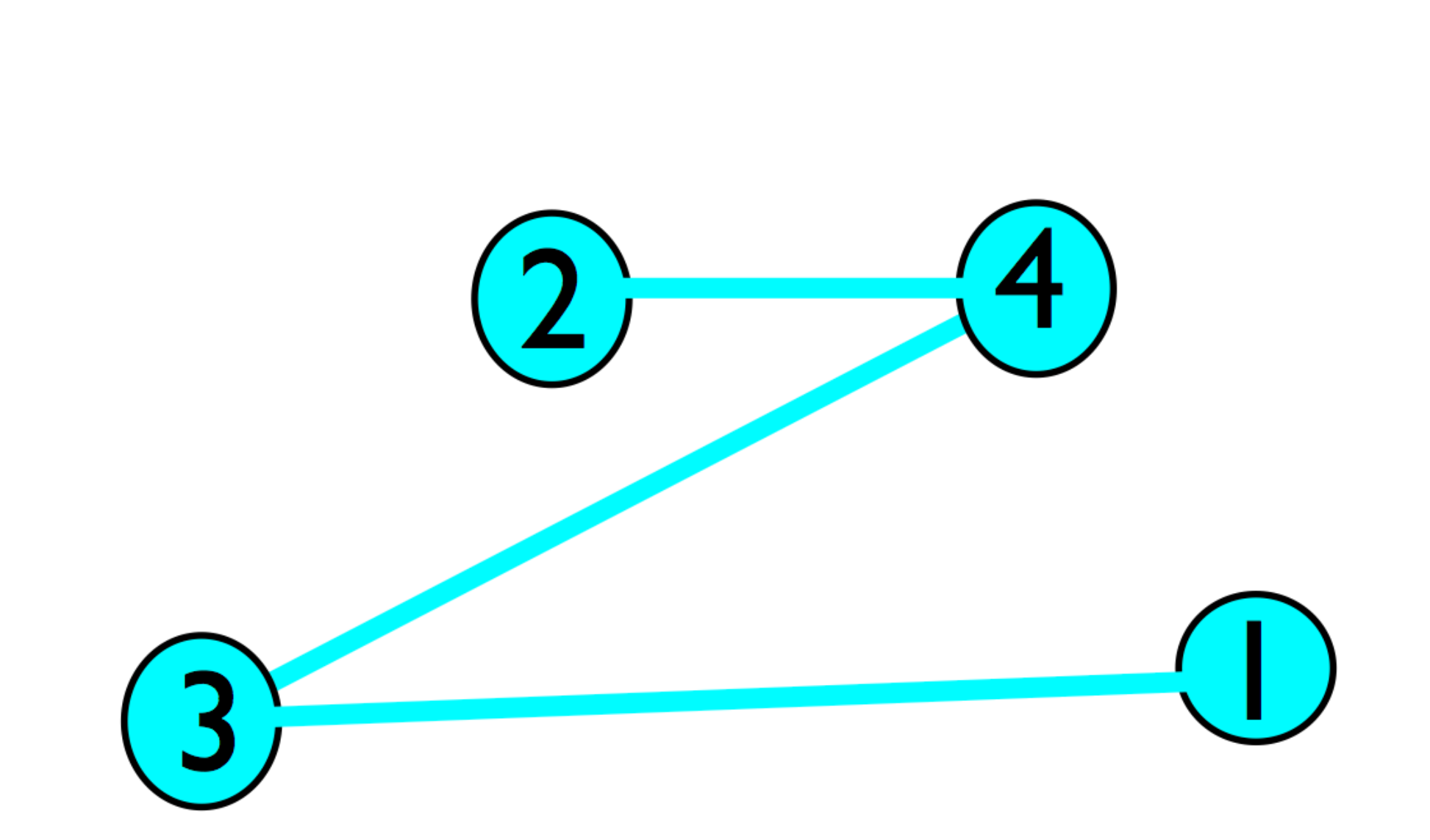}  
\vspace*{-1mm}
\caption{The positive part is connected.}
\label{fig3z}
\end{center}
\vspace{3mm}
\end{figure}

The answer is YES, since, as predicted by the theorem, the configuration 
of positive points and lines is connected (Figure \ref{fig3z}). 

It is natural to wonder how long such a walk can be in the worst case?
Consider for a weighted  configuration of $n$ points
 the graph whose vertices are the positive  points and whose edges are the
pairs of positive points that span a positive line. What is the maximal diameter of such 
a graph?

\medskip

The  theorem contains similar geometric-combinatorial information also for higher ranks.
For instance, fix integers $1<k<r$ and consider a configuration of $n$
points in $\R^{r-2}$. 
A real number weight is assigned to each point, and we assume that
these otherwise generic numbers sum to zero.
By summation each flat of rank $k$ receives a weight, so
 there are  positive and negative $k$-flats.

The theorem guarantees that it is possible to walk in the 
configuration from any positive point 
to any other positive point, never visiting a negative point or moving
across a negative $k$-flat.

The case $k=2$ gives us again the positive point--line graphs,
but now in Euclidean space of arbitrarily high dimension  $\R^{r-2}$.
Again, the diameter question can be asked.

\bigskip

\section{Gorenstein* complexes'}

\subsection{The Gorenstein* property}


The basic topological characterization of Gorenstein* complexes is the
following, derived by Stanley \cite{St77} from the work of Hochster \cite {Ho77}.  

\bthm \label{Thm.Sta}
A simplicial complex $\De$ is {\em Gorenstein*} over {\bf k}  
$$  \Longleftrightarrow \hspace{2mm}  \tH _i (\lk_{\De} (F) ; {\kk})=
\bca \kk, \mbox{ for all $F\in\De$ and $ i =\dim(\lk_{\De}\ F)$}, \\ 
0 \mbox{ for all $F\in\De$ and all $ i < \dim(\lk_{\De}\ F)$}.
\eca
$$
\ethm

For some alternative criteria, see \cite[Theorem  5.1]{St96}. 
The following is a slightly sharper version of criterion $(d)$ of that theorem. 

\bthm Let $\De$ be a simplicial complex, and $\kk=\Z$  \label{Gor1}
 or a field.
Then $\De$ is Gorenstein*  over $\kk$  if and only if 
\ben
\item $\De$ is Cohen-Macaulay over $\kk$, and
\item $\De$ is thin\footnote{A pure complex is said to be {\em thin} if
every face of codimension one lies in exactly two maximal faces.
It is {\em dually connected} if one can walk from any maximal face to 
any other one via steps across codimenson one faces.
It is a {\em pseudomanifold}
if it is thin and dually connected.}
.
\een
If $\kk$ is a field of char $\neq 2$ we also need a third condition: 
\ben
\item[(3)] $\wt{\chi} (\De) =(-1)^{\dim \De}$.
\een
\ethm

Here is a sketch of the proof. The necessity of conditions 
(1), (2) and (3) is easy to deduce from Theorem \ref{Thm.Sta}.
To prove sufficiency, the first observation to be made is that,
in case $\kk=\Z$  \label{Gor11}
 or a field of char $=2$, condition (3)
on the reduced Euler characteristic $\wt{\chi}$ is  
implied by (1) and (2).
 So, for the proof
 of sufficiency condition (3) can be assumed  in all cases.

The second observation is that Cohen-Macaulayness
implies being dually connected , so being Cohen-Macaulay and
thin implies being a pseudomanifold. The two top-dimensional  homology groups
of a pseudo\-manifold  are well known, see e.g. Munkres  \cite{Mu84}.
Using this and verifying for
all links finishes the proof.


\medskip

\subsection{The homotopy Gorenstein* property}

Just like for \cm\ complexes, the Gorenstein* property has a homotopy version.
Somewhat surprisingly, this property turns out to be essentially  
equivalent to being a PL sphere, as shown by Theorem \ref{Gor2a} below.




\bprop   \label{Gor2d}
$\De$ is homotopy Gorenstein* $\Leftrightarrow$ $\De$ is  homotopy Cohen-Macaulay and thin.
\eprop

\bproof 
 We know from Theorem \ref{Gor1} that 
\label{Gor2b}
$$\mbox{ Gorenstein* over $\Z$ \; $\Leftrightarrow$ \, Cohen-Macaulay over 
$\Z$ and  thin.}$$
Adding the condition ``and all links of dimension $\ge 2$ are
simply-connected'' to both sides changes this equivalence to
$$ \mbox{ homotopy Gorenstein* $\Leftrightarrow$ homotopy Cohen-Macaulay 
and  thin.} 
$$
\eproof

 The crucial step in the following proof 
  requires an
 application of the PC-theorem,
 by which we mean the theorem
 of Smale, Freedman and Perelman verifying the generalised
 Poincar\'e conjecture in all dimensions.  Theorem \ref{Gor2c} and
 the usefulness  of  the PC-theorem for its proof
  was discovered a few years ago independently by P. Hersh \cite{He10}
and the author (unpublished).

\bthm   \label{Gor2c}
$\De$ is homotopy Gorenstein* $\Rarr \De$ is homeomorphic to a sphere. 
\ethm

\bproof
The theorem
 is proved by induction on dimension, using that
 $$ \mbox{ if $F\sbseq G \in \De$, then $\lk_{\lk_{\De}(F)}(G\setminus F)=\lk_{\De}(G)$.
}$$
 The claim  is certainly true in low dimensions, for instance if $\dim \De =1$ then
 $\De$ must be a $1$-sphere.
 
 By definition, for $\De$ to be homotopy Gorenstein* means that 
 the link at every face $F\in\De$
  has the homotopy type of a ($\dim \De -1-\dim F)$-sphere. 
 Hence 
 by induction, $\lk_{\De}(F)$ is homeomorphic to a sphere
 of the appropriate dimension 
 for all proper faces $F\in \De$. This shows that $||\De||$ is a 
 topological manifold. Also $\De =\lk_{\De} (\emptyset)$
is by assumption a simply-connected homology sphere.
The PC-theorem therefore implies that $||\De||$ is homeomorphic to a sphere.
\eproof

The inductive procedure of the preceding proof can
be sharpened in many cases, leading to spheres
in the piecewise linear category\footnote{A simplicial
complex which is piecewise linearly homeomorphic to the boundary of a
simplex is called a {\em $PL$ sphere\/}. A {\em combinatorial manifold\/} 
is a triangulation of a topological manifold such that the link at 
every vertex is a $PL$ sphere. A combinatorial manifold which is homeomorphic to a sphere is called a {\em combinatorial sphere}.  PL spheres are combinatorial.

For dimensions $d \neq 4$, a triangulation of the $d$-sphere is
$PL$ if and only if it is combinatorial.
This follows from major work in topology showing that
in these dimensions there is a
unique $PL$ structure for spheres.
For $d=4$ the question whether  
a combinatorial sphere must be 
$PL$  is open. 
}.
The obstruction to being able to say ``in all cases''
is the possibility that non-PL $4$-spheres might turn
up among the links. Defining away this possibility,
we reach this result.

\bthm  \label{Gor2a}
Let $\De$ be a simplicial complex.
Then the following conditions are equivalent:
\ben 
\item $\De$ is homotopy Gorenstein*  and every link of dimension $4$ is a PL sphere
\item $\De$ is a PL sphere. 
\een
\ethm

\bproof
The proof that (1) implies (2) is via the same inductive procedure, now however
keeping track of combinatorial and PL spheres. 
For the inductive step, we have that all proper links in  $\De$ are PL spheres.
That means that $\De$ itseslf is
a combinatorial sphere, and hence PL.

In the other direction,
in a PL  sphere every link of a face 
 is itself a PL sphere, and hence in particular a homotopy sphere.
\eproof

We conclude with a list showing the place of the homotopy Gorenstein*
property in a hierarchy of  near-spherical complexes.
\begin{thm} Let $\De$ be pure. In the following list, \label{Gor3}
each property  implies its successor.
\begin{enumerate}
\item[(a)] $\De$ is thin and shellable 
\item[(b)] $\De$ is a PL sphere


\item[(c)] 
$\De$ is homotopy Gorenstein$^{*}$
\item[(d)] $\De$  is homeomorphic to a sphere
\item[(e)]  
$\De$ is Gorenstein$^{*}$ over ${\Z}$
\item[(f)]  $\De$ is Gorenstein$^{*}$ over ${\bf k}$, for some field {\bf k}

\item[(g)]  $\De$ is a pseudomanifold and $\wt{\chi}(\De) =(-1)^{\dim(\De)}$
\item[(h)] $\De$ is thin
\end{enumerate}
\ethm

Implications in this list that have not been discussed 
in this section are either well known or else 
elementary in view of well known facts. 
All the implications in Theorem \ref{Gor3} are strict, except
possibly (b) $\Rarr$ (c). 






\section{Some personal recollections}

I first met Richard at a conference in Berlin in 1976.
He gave a series of two hour-long talks with title ``Cohen-Macaulay Complexes''.
After that I was hooked. My main interest was and is  in combinatorics, but with strong
side interests in algebra and topology. I was struck by the beauty 
and power of the new area outlined by Richard. 

Then I spent the academic  year
1977--1978 at MIT, which deepened my interest and increased my knowledge. Richard was, as always, very generous with discussions and advice. 
Among a multitude of memories from those days,  let me mention a couple of
minor observations  which seemed  a bit puzzling to me at the time but later
had their explanations. 

One was that Richard's desk seemed to be stuffed 
with packs of sorted index cards.
Often in connection with discussions he would retrieve an index card
and make some annotations. I understood later that this was part of a systematic 
gathering of material for his future books EC1 and EC2, particularly
for the exercises.
At some point I observed that a copy of a Springer Lecture Notes volume 
with the mysterious title ``Toroidal Embeddings'' (\cite{KKMS73}) 
had been lying on his desk for a while.
I asked Richard what the yellow book was all about
and why he was reading it.
At that time I had never heard of toric varieties, 
or of any algebraic-geometric aspect of convex polytopes, so I was
quite amazed by his answer. He said  that he was convinced that
the machinery surrounding these varieties had ingredients that  would one day
add up  to prove the necessity part of the
$g$-conjecture for simplicial polytopes (a 
conjectured characterization of $f$-vectors of simplicial  polytopes).
As we all know, he was right \cite{St80}.

Now,  move fast forward to the spring of 1981. Institut Mittag-Leffler
had a research program then in ``Commutative algebra'',
of which Richard was a participant. During his time in Stockholm Richard
gave a series of eight two-hour lectures on ``Combinatorics and Commutative Algebra'' at Stockholm University, giving a splendid 
overview of the new connections that had been discovered in 
recent years. I had the benefit of taking notes at the lectures and 
writing  up what became the core of Richard' ``green book'' \cite{St83}.
I learned a lot from connecting the dots as
Richard expertly and sometimes rapidly moved from topic to topic.

It so happened that the first Nordic Combinatorics Conference
took place in Norway at the time when Richard was in Stockholm.
The venue was Utstein Abbey, an old 
monastery on an island in the North Sea, just
off the Norwegian coast. We travelled there with a group
of mathematicians by train and boat from Stockholm,
a journey that involved taking night train from Oslo to Stavanger.
The meeting, including frequent refreshment sessions, was 
graciously hosted by Prof. Selmer from Bergen.

For the return trip Richard and I took a different route from the others.
Wanting to experience more of  the majestic Norwegian coast, we took a
local boat north, hopping up the coast to Bergen. After a night there we
continued by train over the mountains  to Oslo, and then on to Stockholm.

There were many hours spent travelling and most of this time we 
discussed mathematics. Somehow the idea to write a book took shape.
 Richard had a lot of material and ideas beyond what
appeared in his Stockholm lectures (or, the first edition of the green book). 
The plan was that Richard
would write about algebraic matters and I would be in charge
of the topological stuff. At that point, I had since a few years 
taken great interest in topological methods in combinatorics and 
poset topology,  an area for which there was much 
nice material but little exposition available at the time. 

As fjords and glaciers passed by outside we worked out more and more detailed plans,
chapter by chapter, leading to the following table of contents.

\vspace{5mm}

{\bf\footnotesize Ring-theoretic and combinatorial aspects of simplicial complexes}

\vspace{3mm}

\hspace*{11mm} {\em\footnotesize Table of contents:}

\vspace{2mm}

\begin{quote}{\em\footnotesize
\begin{enumerate}
\item[1.] Topological preliminaries
\item[2.] Algebraic preliminaries
\item[3.] Cohen-Macaulay complexes
\item[4.] Shellable complexes
\item[5.] $k$-Cohen-Macaulay complexes
\item[6.] Buchsbaum complexes
\item[7.] Gorenstein complexes
\item[8.] Examples
\item[9.] Multigradings
\item[10.] Local cohomology
\item[11.] Canonical modules
\item[12.] Balanced complexes
\item[13.] Bounds
\item[14.] Posets
\item[15.] Topological methods for posets
\item[16.]Group actions
\item[17.] Further developments
\end{enumerate}}
\end{quote}

\vspace{3mm}

As the train finally reached Stockholm we went our separate ways. 
No further work on this book project was ever done.
All that remains is a dusty folder 
with some notes from
our discussions en route from Norway, and
(at least on my side) some fond memories.

\medskip

\small {\bf Acknowledgment:} The pen sketch on the title page is due to artist Berta Hansson. 
It was drawn during the author's thesis defense at 
Stockholm University in May 1979, at which 
Richard Stanley served as opponent

\bibliographystyle{myamsalpha}
\bibliography{rs70}

\end{document}